\documentclass[12pt]{amsart}
\pdfoutput=1
\usepackage[left=1in,right=1in,top=1in,bottom=1in]{geometry}
\usepackage[utf8]{inputenc}
\usepackage{xcolor}
\usepackage[capitalise]{cleveref}
\usepackage{tikz}
\usetikzlibrary{arrows,matrix,positioning}

\usepackage{multirow,tikz, blkarray, mathrsfs,scalerel,stackengine,amssymb}
\usetikzlibrary{arrows,matrix,positioning}
\usepackage{setspace}\setdisplayskipstretch{}

\stackMath
\newcommand\reallywidehat[1]{\savestack{\tmpbox}{\stretchto{\scaleto{ \scalerel*[\widthof{\ensuremath{#1}}]{\kern.1pt\mathchar"0362\kern.1pt}{\rule{0ex}{\textheight}}}{\textheight}}{2.4ex}}\stackon[-6.9pt]{#1}{\tmpbox}}

\newcommand{\C}{\mathbb{C}}
\newcommand{\Q}{\mathbb{Q}}
\newcommand{\Z}{\mathbb{Z}}

\newcommand{\K}{\mathbb{K}}

\newcommand{\reversesproduct}{s_k \dotsm \widehat{s}_{t_r} \dotsm \widehat{s}_{t_1} \dotsm s_1}

\newcommand{\JJ}{\mathcal{J}_{(\lambda \mid \gamma)}}

\newcommand{\MMM}{\mathbb{M}_{\lambdagamma}}

\newcommand{\nudom}{{( \lambda \mid \gamma )}}
\newcommand{\lambdagamma}{{\nudom}}

\newcommand{\mueta}{{(\mu \mid \eta)}}

\newcommand{\mutilde}{{\widetilde{\mu}}}
\newcommand{\lambdatilde}{{\widetilde{\lambda}}}

\newcommand{\lambdatildegamma}{{(\lambdatilde \mid \gamma)}}
\newcommand{\lambdamingamma}{(\lambda^- \mid \gamma)}

\newcommand{\jj}{j_{(\lambda \mid \gamma)}}

\newcommand{\tildepIone}{\widetilde{p}_{I_1}}

\newcommand{\hIone}{h}

\newcommand{\leg}[1]{{\ell_{#1}(\square)}}

\newcommand{\armone}[1]{{\widetilde{a}_{#1}(\square)}}
\newcommand{\armonenu}{{\widetilde{a}_\nu(\square)}}

\newcommand{\armtwo}[1]{{a_{#1}(\square)}}
\newcommand{\armtwonu}{{a_\nu(\square)}}

\newcommand{\symmetricprod}{\prod\limits_{\square \in \lambda^-}}
\newcommand{\nonsymmetricprod}{\prod\limits_{\square \in \gamma}}

\newcommand{\jtwo}{j_{\mathrm{C}}}

\newcommand{\RS}{\mathrm{R_{S}}(x)}
\newcommand{\RNS}{\mathrm{R_{NS}}(x)}
\newcommand{\LL}{\mathrm{LL}(x)}
\newcommand{\RL}{\mathrm{RL}(x)}
\newcommand{\cbox}{b(\square)}

\newcommand{\AAA}{\mathcal{A}^{\lambdagamma}_{\mueta}}

\newcommand{\Atq}{\mathbb{A}_{t,q}}
\newcommand{\At}{\mathbb{A}_t}
\newcommand{\id}{\mathsf{e}}

\newcommand{\PFH}{\mathrm{PFH}}

\newcommand{\Hilb}{\mathrm{Hilb}}

\newcommand{\Hom}{\mathrm{Hom}}

\newcommand{\T}{\mathbb{T}}

\newcommand{\wt}{\mathrm{wt}}

\newcommand{\comment}[1]{}

\newcommand{\Y}{\mathbb{Y}}
\newcommand{\dg}{\mathrm{dg}}
\newcommand{\dga}{\widehat{\dg}}

\newcommand{\la}{\lambda}
\newcommand{\ga}{\gamma}

\newcommand{\tH}{\tilde{H}}

\newcommand{\midd}{\,|\,}

\newcommand{\dd}{\mathsf{d}}
\newcommand{\TT}{\mathsf{T}}
\newcommand{\yy}{\mathsf{y}}
\newcommand{\zz}{\mathsf{z}}

\newcommand{\cJ}{\mathcal{J}}

\newcommand{\armm}{\mathrm{arm}}
\newcommand{\legg}{\mathrm{leg}}

\theoremstyle{plain}
\newtheorem{thm}{Theorem}[section]
\newtheorem{cor}[thm]{Corollary}
\newtheorem{lem}[thm]{Lemma}
\newtheorem{prop}[thm]{Proposition}
\newtheorem{conj}[thm]{Conjecture}
\newtheorem*{conj*}{Conjecture}

\theoremstyle{definition}
\newtheorem{dfn}[thm]{Definition}
\newtheorem{ex}[thm]{Example}
\newtheorem{rem}[thm]{Remark}

\numberwithin{equation}{section}

\title{A geometric realization of partially-symmetric Macdonald polynomials}
\author{Ben Goodberry}
\address{B.G.: Department of Mathematical Science, Henson Science Hall, Salisbury University, 1101 Camden Avenue, Salisbury, MD 21801}
\email{bngoodberry@salisbury.edu}
\author{Daniel Orr}
\address{D.O.: Department of Mathematics (0123),
460 McBryde Hall, Virginia Tech,
225 Stanger Street,
Blacksburg, VA 24061-1026}
\email{dorr@math.vt.edu}
\date{\today}

\begin{document}

\maketitle

\begin{abstract}

We formulate a precise conjecture relating integral form partially-symmetric Macdonald polynomials and the parabolic flag Hilbert schemes of Carlsson, Gorsky, and Mellit. This extends, in an explicit fashion, Haiman's realization of modified Macdonald symmetric functions via Hilbert schemes of points in the plane. As evidence for our conjecture we prove that it is compatible with the action of certain elements in Carlsson and Mellit's algebra $\Atq$, including degree $1$ Pieri formulas.

\end{abstract}

\section{Introduction}

Let $\Lambda$ be the algebra of symmetric functions in $x_1,x_2,\dotsc$ over $\K=\Q(q,t)$. A key structure in Haiman's geometric realization of Macdonald symmetric functions \cite{H1} (see also \cite{H2}) is a vector space isomorphism
\begin{align}\label{E:Phi-small}
\bigoplus_{n\ge 0} K_\T(\Hilb_n)_\mathrm{loc} \overset{\Phi_0}{\longrightarrow} \Lambda
\end{align}
between $\Lambda$ and the localized equivariant $K$-groups of the Hilbert schemes $\Hilb_n$ of points in the plane $\C^2$, where $\T$ is the natural two-dimensional torus acting on the plane. Under this isomorphism, the modified Macdonald symmetric functions $\tH_\mu$ are realized as the images $\Phi_0([I_\mu])$ of the $K$-theory classes corresponding to torus-fixed points $I_\mu$ in $\sqcup_{n\ge 0} \Hilb_n$. It is of great significance that the map $\Phi_0$ has a genuinely geometric origin---namely, the Procesi bundle---which exhibits Macdonald positivity.

Work of Schiffmann and Vasserot \cite{SV} offers another point of view on the map $\Phi_0$: it is an intertwining map between two irreducible representation of the elliptic Hall algebra $\mathcal{E}$. While this point of view is not strong enough to establish Macdonald positivity, it does provide a way to uniquely characterize the map $\Phi_0$, up to a scalar.

By methods akin to those of \cite{SV}, Carlsson, Gorsky, and Mellit \cite{CGM} costruct an extension of the map $\Phi_0$, with Carlsson and Mellit's algebra $\Atq$ \cite{CM} in place of the elliptic Hall algebra $\mathcal{E}$. The extension $\Phi$ maps between the following larger spaces
\begin{align}\label{E:Phi-big}
\bigoplus_{n\ge 0}\bigoplus_{0\le k\le n} K_\T(\PFH_{n,n-k})_\mathrm{loc} \overset{\Phi}{\longrightarrow} \bigoplus_{k\ge 0}\Lambda\otimes\K[y_1,\dotsc,y_k]
\end{align}
both of which afford an irreducible action of the algebra $\Atq$ with a canonical generator. Here $\PFH_{n,n-k}$ is the parabolic flag Hilbert scheme, also introduced in \cite{CGM},  which parametrizes flags of ideals $I_n\subset I_{n-1} \subset \dotsm \subset I_{n-k}$ in $\C[x,y]$ with the property that $\dim \C[x,y]/I_s = s$ and $yI_{n-k}\subset I_n$. The action of $\Atq$ on the direct sum of $K$-groups of these spaces is constructed in \cite{CGM} by geometrically-defined operators, while the right-hand side of \eqref{E:Phi-big} affords the polynomial representation of \cite{CM}. Of course, the requirement that $\Phi$ intertwine these two actions is sufficient to uniquely characterize it, up to a scalar.

Our work is motivated by the following question: \textit{what are the images of torus fixed-point classes in $\PFH_{n,n-k}$ under the isomorphism $\Phi$?} One knows that $\Phi$ restricts to \eqref{E:Phi-small} on the $k=0$ summands (up to a minor twist by a sign and a monomial in $q,t$), where the torus fixed point-classes are sent to modified Macdonald symmetric functions $\tH_\mu$. We propose an explicit answer for all $k$ involving a new family of modified partially-symmetric Macdonald functions $\tH_{(\la|\ga)}$.

\begin{conj*}[see Conjecture \ref{C:main}]
Under the $\Atq$-module isomorphism 
$$ \bigoplus_{n\ge 0}\bigoplus_{0\le k\le n} K_\T(\PFH_{n,n-k})_{\mathrm{loc}} \overset{\Phi}{\longrightarrow} \bigoplus_{k\ge 0}\Lambda\otimes\K[y_1,\dotsc,y_k], $$
the fixed-point classes are sent to modified partially-symmetric Macdonald polynomials, as follows:
\begin{align*}
H_{\mu,w} \mapsto t^{n(\mathrm{sort}(\lambda,\gamma))+|(\lambda|\gamma)|}\mathcal{J}_{(\lambda|\gamma)}^{w_0}\Big(\frac{X}{t^{-1}-1}\,\big|\,  y\Big)^*=:\tH_{(\la|\ga)}.
\end{align*}
\end{conj*}

For the sake of orientation, let us describe a few pieces of the notation used above, leaving the rest for the main body of the paper. The fixed-points $I_{\mu,w}\in\PFH_{n,n-k}$ are indexed by pairs $(\mu,w)$ such that $\mu$ is a partition of size $n$ and $w$ is an ordered horizontal strip of size $k$ in $\mu$; the classes $H_{\mu,w}$ are explicit scalar multiples of the $[I_{\mu,w}]$. The partially-symmetric integral form Macdonald functions $\cJ_{(\la|\ga)}^{w_0}(X\midd  y_1,\dotsc,y_k)\in\Lambda\otimes\K[y_1,\dotsc,y_k]$ are defined in \eqref{E:J-w0} by partially-symmetrizing the (stable) nonsymmetric Macdonald polynomials, multiplying by explicit scalars introduced in \cite{G} and \cite{L}, and then performing a twist by the long element of $S_k$. The $\cJ_{(\la|\ga)}^{w_0}(X\midd y)$---and their cousins $\tH_{(\la|\ga)}$---are indexed by pairs $\lambdagamma$ consisting of a partition $\lambda$ and a composition $\gamma\in(\Z_{\ge 0})^k$ such that the total size $|\lambda|+|\gamma|$ is $n-k$. An explicit bijection $\phi$ between these indexing sets, so that $\phi(\mu,w)=\lambdagamma$, is recalled in Section~\ref{S:Bij} following \cite{CGM}. Finally, $X\mapsto X/(1-t)$ denotes the familiar plethystic transformation on $\Lambda$ from Macdonald theory, and $f(X\midd y)^*=f(X\midd y)|_{t\mapsto t^{-1}}$.

Thus our conjecture closely resembles---and clearly extends---the construction of modified Macdonald functions $\tH_\mu$ from the usual integral form Macdonald functions.

In this paper, we make substantial progress toward a proof of Conjecture~\ref{C:main} by showing that it is compatible with the action of certain $\Atq$-generators; see Theorem~\ref{T:main}. We will complete the proof in forthcoming work with M. Bechtloff Weising by means of the limit Cherednik operators \cite{IW,BW}.

It is important to mention that the map $\Phi$ of \cite{CGM}, while uniquely characterized as an $\Atq$-intertwiner, has no genuinely geometric construction at present. In light of our conjecture and Lapointe's positivity conjectures \cite{L} for partially-symmetric Macdonald polynomials, it natural to hope for a geometric realization of $\Phi$ exhibiting positivity.

\section*{Acknowledgements}

We thank Milo Bechtloff Weising, Nicolle Gonz\'{a}lez, Mark Haiman, Bogdan Ion, Luc Lapointe, Jennifer Morse, and Mark Shimozono for helpful conversations. D.O. was partially suppported by the Simons Foundation (Award 638577).

\section{Partially symmetric Macdonald polynomials}

\subsection{Diagrams}

In constructions related to compositions and their diagrams, we mostly follow the conventions of \cite{HHL} (except that we write $\dg$ instead of $\dg'$ for their diagrams). Compositions are tuples $\nu=(\nu_1,\dotsc,\nu_n)\in(\Z_{\ge 0})^n$. A composition $\nu$ is a partition (with at most $n$ parts) if its entries are weakly decreasing. Let $\Y_n\subset (\Z_{\ge 0})^n$ be the set of partitions with at most $n$ parts. For any composition $\nu\in(\Z_{\ge 0})^n$, let $\mathrm{sort}(\nu)$ be the unique partition in $\Y_n$ obtained by rearranging the parts of $\nu$. For $\nu\in(\Z_{\ge 0})^n$, let $|\nu|=\nu_1+\dotsm+\nu_n$. 

The diagram of a composition $\nu\in(\Z_{\ge 0})^n$ is the subset $\dg(\nu)\subset (\Z_{\ge 0})^2$ given by
$$ \dg(\nu)=\{(i,j) \mid 1\le i\le n, 1\le j\le\nu_i\}. $$
We view the parts $\nu_i$ of $\nu$ as columns in $\dg(\nu)$. Elements of $\dg(\nu)$ are called boxes of $\nu$. We identify a composition with its diagram and write $(i,j)\in\nu$ to mean $(i,j)\in \dg(\nu)$. We also define the following subsets of $\dg(\nu)$:
\begin{align*} 
\dg_r(\nu) &= \{(i,j) \in \nu \mid j=r \}
\end{align*}

The leg and arm lengths of a box $\square=(i,j)\in\nu$ are defined by
\begin{align*}
\leg{\nu} &= \nu_i -j,\\
\armtwonu &= \# \{ 1 \leq r < i \mid j\leq \nu_r \leq \nu_i\} + \# \{ i<r \leq n \mid j-1\leq \nu_r < \nu_i\}.
\end{align*}
We will also use the following alternate versions of arm length (see Example~\ref{DiagramExample}):
\begin{align*}
\armonenu &= \# \{1\leq r < i \mid j \leq \nu_r \leq \nu_i \} + \# \{i<r \leq n \mid j \leq \nu_r < \nu_i \}.
\end{align*}

Let $\Y=\cup_{n\ge 0} \Y_n$ be the set of all partitions. For any $\nu\in\Y$, denote by $\nu'\in\Y$ the transposed partition and let $\dg'(\nu)=\dg(\nu')$. Thus $\dg'(\nu)$ is the usual (French) diagram of partition $\nu$ given by rows of boxes of lengths $\nu_1,\nu_2,\dotsc$ in the first quadrant. Let $\mathrm{leg}_\nu$ and $\mathrm{arm}_\nu$ denote the usual leg and arm functions on $\dg'(\nu)$ as defined in, e.g., \cite{H1}. (For a box $\square\in\dg'(\nu)$, we have $\mathrm{leg}_\nu(\square)=\ell_{\nu'}(\square)$, and $\mathrm{arm}_\nu(\square)$ is related to the function $a_{\nu_-}$ for the weakly increasing rearrangement $\nu_-$ of $\nu$.)

For any $\nu\in\Y$, define $n(\nu)=\sum_i (i-1)\nu_i$.



\subsection{Split diagrams}

We write $\nu = (\lambda|\gamma)$ to indicate a splitting of a composition $\nu$ into two compositions $\lambda$ and $\gamma$, i.e., $\nu=(\nu_1,\dotsc,\nu_n)=(\lambda_1,\dotsc,\lambda_m,\gamma_1,\dotsc,\gamma_k)$ where $\lambda=(\lambda_1,\dotsc,\lambda_m)\in(\Z_{\ge 0})^m$ and $\gamma=(\gamma_1,\dotsc,\gamma_k)\in(\Z_{\ge 0})^k$ and $n=m+k$. In this situation, we will often assume in addition that $\lambda$ is a partition or the reverse of a partition.

Starting with $\nu=(\lambda | \gamma)$ as above and $\lambda$ a partition, we construct an augmented diagram $\dga(\nu)$ as follows. First, we form $\nu^- := (\lambda^- | \gamma)$, where $\lambda^-$ is the weakly increasing rearrangement of $\lambda$. The augmented diagram associated with $\nu$ is then defined as follows
$$ \dga(\nu)=\dg(\nu^-)\cup \{(m+1,0),\dotsc,(n,0)\}. $$ 
Note that the additional boxes $(m+1,0),\dotsc,(n,0)$ lie directly below the columns of $\dg(\nu^-)$ corresponding to $\gamma$. 

We call the subsets of $\dg(\nu^-)$ corresponding to $\lambda^-$ and $\gamma$ the symmetric and nonsymmetric parts of the diagram, respectively. We will use different arm functions for boxes in these two parts of $\dg(\nu^-)$, as follows. For boxes in the symmetric part of $\dg(\nu^-)$, we will use $\widetilde{a}_{\nu^-}$. For boxes in the nonsymmetric part, we will use the arm function $a_{\nu^-}$. These arm functions have interpretations as counting certain boxes in $\dga(\nu^-)$, as we illustrate in the example below. 




\begin{ex}\label{DiagramExample}

Let $\nu = (3,1 \mid 2, 1, 3, 0, 1)$. Then $\nu^-=(1,3\mid 2,1,3,0,1)$ and the augmented diagram $\dga(\nu^-)$ is
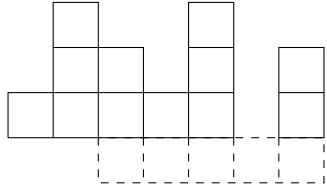
\begin{center}
\begin{tikzpicture}[scale=.6]
\draw (0,0) -- (0,1) -- (5,1) -- (5,0) -- (0,0);
\draw (1,0) -- (1,3) -- (2,3) -- (2,0);
\draw (1,2) -- (3,2) -- (3,0);
\draw (4,0) -- (4,3) -- (5,3) -- (5,1);
\draw (4,2) -- (5,2);
\draw (6,0) -- (6,2) -- (7,2) -- (7,0) -- (6,0);
\draw (6,1) -- (7,1);
\draw[dashed] (2,0) -- (2,-1) -- (7,-1) -- (7,0) -- (2,0);
\draw[dashed] (3,0) -- (3,-1);
\draw[dashed] (4,0) -- (4,-1);
\draw[dashed] (5,0) -- (5,-1);
\draw[dashed] (6,0) -- (6,-1);
\end{tikzpicture}
\end{center}
where the boxes added to $\dg(\nu^-)$ are dashed.

Let us illustrate the boxes that are counted for each of the arms and legs. For $u=(3,1)$, a box in the nonsymmetric part of $\dga(\nu^-)$, the following boxes marked $\ell$ and $a$ contribute to $\ell_{\nu^-}(u)=1$ and $a_{\nu^-}(u)=1+2=3$:
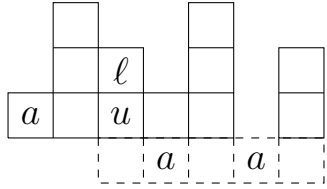
\begin{center}
\begin{tikzpicture}[scale=.6]
\draw (0,0) -- (0,1) -- (5,1) -- (5,0) -- (0,0);
\draw (1,0) -- (1,3) -- (2,3) -- (2,0);
\draw (1,2) -- (3,2) -- (3,0);
\draw (4,0) -- (4,3) -- (5,3) -- (5,1);
\draw (4,2) -- (5,2);
\draw (6,0) -- (6,2) -- (7,2) -- (7,0) -- (6,0);
\draw (6,1) -- (7,1);
\draw[dashed] (2,0) -- (2,-1) -- (7,-1) -- (7,0) -- (2,0);
\draw[dashed] (3,0) -- (3,-1);
\draw[dashed] (4,0) -- (4,-1);
\draw[dashed] (5,0) -- (5,-1);
\draw[dashed] (6,0) -- (6,-1);
\node at (.5,.5) {\large $ a$};
\node at (2.5,.5) {\large $ u$};
\node at (2.5,1.5) {\large $ \ell$};
\node at (3.5,-.5) {\large $ a$};
\node at (5.5,-.5) {\large $ a$};
\node at (6.5,-.5) {};
\end{tikzpicture}
\end{center}
For $v = (2,1)$, in the symmetric part of the diagram, the boxes that contribute to $\ell_{\nu^-}(u)=2$ and $\widetilde{a}_{\nu^-}(v)=1+3=4$ are:
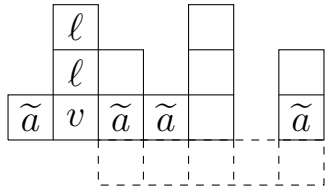
\begin{center}
\begin{tikzpicture}[scale=.6]
\draw (0,0) -- (0,1) -- (5,1) -- (5,0) -- (0,0);
\draw (1,0) -- (1,3) -- (2,3) -- (2,0);
\draw (1,2) -- (3,2) -- (3,0);
\draw (4,0) -- (4,3) -- (5,3) -- (5,1);
\draw (4,2) -- (5,2);
\draw (6,0) -- (6,2) -- (7,2) -- (7,0) -- (6,0);
\draw (6,1) -- (7,1);
\draw[dashed] (2,0) -- (2,-1) -- (7,-1) -- (7,0) -- (2,0);
\draw[dashed] (3,0) -- (3,-1);
\draw[dashed] (4,0) -- (4,-1);
\draw[dashed] (5,0) -- (5,-1);
\draw[dashed] (6,0) -- (6,-1);
\node at (1.5,.5) {\large $v$};
\node at (1.5,1.5) {\large $\ell$};
\node at (1.5,2.5) {\large $\ell$};
\node at (.5,.5) {\large $\widetilde{a}$};
\node at (2.5,.5) {\large $\widetilde{a}$};
\node at (3.5,.5) {\large $\widetilde{a}$};
\node at (6.5,.5) {\large $\widetilde{a}$};
\end{tikzpicture}
\end{center}
\end{ex}

\subsection{Nonsymmetric Macdonald polynomials}

Let $\K=\Q(q,t)$ and $n>0$. For any composition $\nu\in(\Z_{\ge 0})^n$, denote by $E_\nu =E_\nu(x;q,t)\in \K[x_1,\dotsc,x_n]$ the nonsymmetric Macdonald polynomial of type $GL_n$ as defined in  \cite{HHL}.

Let the symmetric group $S_n$ act in the natural way on $\Z^n$ and on $\K[x_1,\dotsc,x_n]$. For $1\le i<n$, denote by $s_i$ the simple transposition $(i,i+1)\in S_n$. Define the Demazure-Lusztig operators, acting on $\K[x_1,\dotsc,x_n]$, by
\begin{align}\label{E:DL}
T_i &= ts_i + \frac{(t-1)x_{i+1}}{x_{i+1}-x_i}(1-s_i)=\frac{x_{i+1}-tx_i}{x_{i+1}-x_i}s_i + \frac{(t-1)x_{i+1}}{x_{i+1}-x_i}, \qquad 1\le i<n.
\end{align}
These satisfy the same braid relations as the $s_i$ and the following quadratic relations:
\begin{equation}\label{quadratic}
(T_i-t)(T_i+1)=0.
\end{equation}
For any reduced expression $w=s_{i_1}\dotsm s_{i_\ell}\in S_n$, we unambiguously define the operator $T_w=T_{i_1}\dotsm T_{i_\ell}$. 


\begin{lem}[{\cite[(17)]{HHL}}]
If $\nu_{i}>\nu_{i+1}$ for some $1\le i<n$, then
\begin{align}
\label{E:T-on-E-up}
T_i E_\nu = E_{s_i(\nu)}-\frac{1-t}{1-q^{\ell_\nu(\square)+1}t^{a_\nu(\square)}}E_\nu
\end{align}
where $\square=(i,\nu_{i+1}+1)\in\dg(\nu)$.
\end{lem}

\subsection{Partially symmetric Macdonald polynomials}

Suppose $m,k\ge 0$ and $n=m+k>0$. In this situation, we regard $S_m$ as the subgroup of $S_n$ fixing the elements $\{m+1,\dotsc,n\}$. Define the partial Hecke symmetrizer
$$
P^+_m=\sum_{w\in S_m} T_w.
$$
For $\lambda\in(\Z_{\ge 0})^m$, let $(S_m)_\lambda\subset S_m$ be its stabilizer and
\begin{align*}
S_\lambda(t)=\sum_{w\in (S_m)_\lambda} t^{\ell(w)}.
\end{align*}

\begin{dfn}
For a partition $\lambda\in\Y_m$ and a composition $\gamma\in (\Z_{\ge 0})^k$, the partially symmetric Macdonald polynomial $P_{(\lambda|\gamma)}=P_{(\lambda|\gamma)}(x;q,t)$ is defined by
$$
P_{(\lambda|\gamma)} = \frac{P^+_m \cdot E_{(\lambda|\gamma)}}{S_\lambda(t)}.
$$
\end{dfn}

The $P_{(\la|\ga)}$ form a basis for the space $\K[x_1,\dotsc,x_n]^{S_m}$ of partially-symmetric polynomials.

\begin{dfn}
The integral form partially symmmetric Macdonald polynomial $J_{(\lambda|\gamma)}=J_{(\lambda|\gamma)}(x;q,t)$ is the scalar multiple
$$
\cJ_{(\lambda|\gamma)}= j_{(\lambda|\gamma)}P_{(\lambda|\gamma)}
$$
where
$$ 
\jj := \symmetricprod (1-q^{\ell(\square)}t^{\tilde{a}(\square) + 1}) \nonsymmetricprod(1-q^{\ell(\square) + 1}t^{a(\square) + 1})
$$
and all arms and legs are taken in the augmented diagram $\dga(\nu^-)$ where $\nu^-=(\lambda^-\mid \gamma)$.
\end{dfn}

It is shown in \cite{G} (and asserted in \cite{L}) that $\cJ_{(\la|\ga)}\in\Z[q,t][x_1,\dotsc,x_n]$.




For later use, we record the action of $T_{m+i}$ for $i=1,\dotsc,k$ on the integral forms:
\begin{align}\notag
T_{m+i}\cJ_{(\lambda|\gamma)} &= \frac{j_{(\lambda|\gamma)}}{j_{(\lambda|s_i(\gamma))}}\cJ_{(\lambda|s_i(\gamma))}-\frac{1-t}{1-q^{l(u)+1}t^{a(u)}}\cJ_{(\lambda|\gamma)}\\
&= \frac{1-q^{l(u)+1}t^{a(u)+1}}{1-q^{l(u)+1}t^{a(u)}}\cJ_{(\lambda|s_i(\gamma))}-\frac{1-t}{1-q^{l(u)+1}t^{a(u)}}\cJ_{(\lambda|\gamma)}.\label{E:TJ}
\end{align}
where arms and legs are taken in $\dg(\la|\ga)$ such that $\ga_i>\ga_{i+1}$ and $u=(m+i,\gamma_{i+1}+1)$.

\subsection{Stability}

For any $n=m+k$ with $m,k\ge 0$, we make the identification
\begin{align}\label{E:PSI}
\K[x_1,\dotsc,x_n]^{S_m}\cong \K[x_1,\dotsc,x_m]^{S_m}\otimes \K[y_1,\dotsc,y_k]
\end{align}
by mapping $x_1,\dotsc,x_m$ to themselves and $x_{m+1},\dotsc,x_n$ to $y_1,\dotsc,y_k$, respectively. We will write $f(x\mid y)$ for the image of $f\in\K[x_1,\dotsc,x_n]^{S_m}$ under this identification.

With this notation in place, we can formulate the stability of partially symmetric Macdonald polynomials, which is proved in \cite{L} and \cite{G}.

\begin{prop}
For any $\la\in\Y_m$ and $\ga\in(\Z_{\ge 0})^k$, we have
$$ P_{(\la, 0 \, \mid\, \ga)}(x,0\mid y) = P_{\lambdagamma}(x \midd y).$$
\end{prop}

Since $j_{(\la,0\, \mid \, \ga)}=j_{\lambdagamma}$, the integral forms $\cJ_{\lambdagamma}$ have the same stability property. Taking $m\to\infty$, we write $P_{\lambdagamma}=P_{\lambdagamma}(X\midd y)$ and $\cJ_{\lambdagamma}=\cJ_{\lambdagamma}(X\midd y)$ for the corresponding elements of $\Lambda\otimes\K[y_1,\dotsc,y_k]$, where $\Lambda$ is the $\K$-algebra of symmetric functions in $x_1,x_2,\dotsc$. As $\la$ ranges over all partitions and $\ga\in(\Z_{\ge 0})^k$, these form bases of $\Lambda\otimes\K[y_1,\dotsc,y_k]$.

By slight abuse of notation, we now denote by $T_i$ for $i=1,\dotsc,k-1$ the Demazure-Luzstig operators \eqref{E:DL} acting on the variables $y_1,\dotsc,y_k$ in $\Lambda\otimes\K[y_1,\dotsc,y_k]$. Under the identification \eqref{E:PSI}, this $T_i$ corresponds to $T_{m+i}$ in the variables $x_{m+i},\dotsc,x_n$.

We define the following variant of the integral form Macdonald functions in $\Lambda\otimes\K[y_1,\dotsc,y_k]$:
\begin{align}\label{E:J-w0}
\cJ^{w_0}_{\lambdagamma}(X\midd y) = t^{-\ell(w_0)} w_0 T_{w_0} \cJ_{\lambdagamma}(X\midd y)
\end{align}
where $w_0\in S_k$ is the long element.

\section{Carlsson-Mellit algebra}

All material from this section is due to \cite{CM} and \cite{CGM}.

\subsection{The algebra}
Let $\Atq$ be the Carlsson-Mellit algebra as defined in \cite{CGM}; note that compared to \cite{CGM} we swap the roles of $q$ and $t$. By definition, $\Atq$ is an associative $\K$-algebra generated by orthogonal idempotents $\id_0,\id_1,\dotsc$ and elements $$\dd_-,\dd_+,\dd_+^*,\TT_1,\TT_2,\dotsc,\yy_1,\yy_2,\dotsc,\zz_1,\zz_2,\dotsc$$ subject to certain relations. We only state two of the relations in order to make our difference in conventions clear: $\zz_1\dd_+\id_k=-qt^{k+1}\yy_1\dd_+^*\id_k$ and $(\TT_i-1)(\TT_i+t)=0$ for all $i$. 

Let $\At\subset\Atq$ be the subalgebra generated by the elements $\dd_-,\dd_+,\TT_1,\TT_2,\dotsc,\yy_1,\yy_2,\dotsc$ and the orthogonal idempotents $\id_0,\id_1,\dotsc$. (In fact, the $\yy_1,\yy_2,\dotsm$ are expressible in terms of the other generators.)

\subsection{Polynomial representation}

The polynomial representation of $\Atq$ was introduced in \cite{CM}. We use the conventions of \cite{CGM}, where it was shown that the polynomial representation is isomorphic 
to the geometric representation described in the next subsection.

Recall that $\Lambda$ is the $\K$-algebra of symmetric functions in $x_1,x_2,\dotsc$. The polynomial representation space is $V=\oplus_{k\ge 0} V_k$ where $V_k=\Lambda\otimes\K[y_1,\dotsc,y_k]$. The idempotent $\id_k$ acts by the natural projection $\id_k : V \to V_k$. On $V_k=\id_k V$, the operators $\TT_1,\dotsc,\TT_{k-1}$ and $\yy_1,\dotsc,\yy_k$ are given by
\begin{align*}
\TT_i \cdot f(X\midd y) &= \frac{(t-1)y_{i}}{y_{i+1}-y_i}f+\frac{y_{i+1}-ty_i}{y_{i+1}-y_i}s_if\\
\yy_i \cdot f &= y_i f
\end{align*}
where $f=f(X\midd y)$ and $s_i$ permutes $y_i$ and $y_{i+1}$.
Note that $\TT_i$ are closely related to (but not the same as) the Demazure-Lusztig operators $T_i$; we will make the relation between these operators precise in Section~\ref{S:T} below.

The operators $\dd_- : V_{k+1}\to V_k$ and $\dd_+: V_k \to V_{k+1}$ are given by
\begin{align*}
\dd_-\cdot f &=-[y_{k+1}^{-1}]\Big(f(X-(t-1)y_{k+1})\Omega(-y_{k+1}^{-1}X)\Big)\\ 
\dd_+\cdot F &= \TT_1\dotsm \TT_k\cdot f(X+(t-1)y_{k+1}),
\end{align*}
where $[y^s]A(y)=A_s$ for a formal series $\sum_{s\in\Z} A_s y^s$, and $\Omega(-yX)=\sum_{i\ge 0} (-y)^i e_i(X)$ in terms of the elementary symmetric functions $e_i$. We use plethystic notation for the $\K[y_1,\dotsc,y_k]$-linear automorphism $f\mapsto f(X\pm (t-1)y_{k+1})$ of $V_k$, which is determined by its values on power sums $p_i$:
\begin{align*}
p_i(X\pm(a-b)) = p_i(X)\pm (a^i-b^i).
\end{align*}
These formulas define $V$ as an $\At$-module; in fact, $V\cong \At \id_0$ as $\At$-modules by \cite{CM}. For the action of the remaining generators, we refer the reader to \cite{CGM}.

One easily checks that the operator of multiplication by the symmetric function $e_1(X)$ on $f=f(X\midd y)\in V_k$ is expressed in terms of the generators of $\At$ as follows:
\begin{align*}
e_1(X) f = \dd_- \TT_k^{-1}\dotsm \TT_1^{-1} \dd_+\cdot f.
\end{align*}

\begin{rem}
We have corrected a minor typo in the action of $\TT_i$ on $V$ made in \cite{CGM}.
\end{rem}


\subsection{Parabolic flag Hilbert schemes}\label{S:PFH}

For any integers $k$ and $n$ such that $0\le k\le n$, let $\PFH_{n,n-k}$ denote the parabolic flag Hilbert scheme of points in $\C^2$ introduced in \cite{CGM}. The points of $\PFH_{n,n-k}$ are flags of ideals
\begin{align*}
I_n \subset I_{n-1} \subset \dotsm \subset I_{n-k} \subset \C[x,y],
\end{align*}
such that $I_{s}$ has codimension $s$ in $\C[x,y]$ for all $s$, and
\begin{align*}
yI_{n-k}\subset I_n.
\end{align*}
The space $\PFH_{n,n-k}$ has the structure of a smooth algebraic variety of dimension $2n-k$. We note that $\PFH_{n,n}$ is the usual Hilbert scheme $\Hilb_n$ of $n$ points in $\C^2$, while $\PFH_{n,0}$ is isomorphic to $\C^n$.

The group $\T=\C^*\times\C^*$ acts on $\C[x,y]$, $\Hilb_n$, and $\PFH_{n,n-k}$ via its natural action on $\C^2$. Let $q,t\in\Hom(\T,\C^*)$ be the $\T$-weights of $y$ and $x$, respectively; a monomial $y^r x^c$ has $\T$-weight $q^r t^c$. We visualize $y^r x^c$ as being located in row $r$ and column $c$ of an infinite array of unit boxes in the first quadrant, with indices beginning at $0$. We also refer to $q^rt^c$ as the $\T$-weight of the box $\square=(c,r)$, denoted $\wt(\square)=q^rt^c$.

The $\T$-fixed points of $\Hilb_n$ are precisely the monomial ideals
\begin{align*}
I_\mu=(x^{\mu_1}y^0,x^{\mu_2}y^1,\dotsc),
\end{align*}
where $\mu=(\mu_1,\mu_2,\dotsc)\in\Y$ is a partition of size $|\mu|=n$. The $\T$-fixed points of $\PFH_{n,n-k}$ are given by flags of monomial ideals
\begin{align*}
I_{\mu^{(n)}} \subset I_{\mu^{(n-1)}} \subset \dotsm \subset I_{\mu^{(n-k)}}\subset \C[x,y] 
\end{align*}
where $\mu^{(s)}$ is a partition of size $s$ for each $s$ and $y I_{\mu^{(n-k)}}\subset I_{\mu^{(n)}}$. In order to satisfy these conditions, the indexing partitions must satisfy
\begin{align}\label{E:chain-partitions}
\dg'(\mu^{(n)}) \supset \dg'(\mu^{(n-1)}) \supset \dotsm \supset \dg'(\mu^{(n-k)}),
\end{align}
with $\square_s=\dg'(\mu^{(s)})\setminus\dg'(\mu^{(s-1)})$ having size $1$ for each $s$, and with $\dg'(\mu^{(n)})\setminus\dg'(\mu^{(n-k)})$ being a horizontal strip (i.e., having at most one box in each column.)

Following \cite{CGM}, we will denote by $I_{\mu,w}$ the $\T$-fixed point of $\PFH_{n,n-k}$ determined by a chain of partitions \eqref{E:chain-partitions}, where $\mu=\mu^{(n)}$ and $w=(w_1,\dotsc,w_k)$ with $w_t=\wt(\square_{n-t+1})$ for $1\le t\le k$. Thus $\mu$ is the biggest partition in the chain, and the $w_t$ are the weights of the boxes whose removal determines the rest of the flag, in the order of removal.

\begin{ex}
The $\T$-fixed point $I_{(3,1),(t^2,q,t)}$ in $\PFH_{4,1}$ is given by the flag
$I_{(3,1)}\subset I_{(2,1)} \subset I_{(2)} \subset I_{(1)}.$ 

\begin{center}
    \begin{tikzpicture}[scale=0.7,fill opacity=.3,every node/.style={scale=0.7}]
\draw (0,0) -- (3,0);
\draw (0,1) -- (3,1);
\draw (0,2) -- (1,2);
\draw (0,0) -- (0,2);
\draw (1,0) -- (1,2);
\draw (2,0) -- (2,1);
\draw (3,0) -- (3,1);
\node[fill opacity = 1] at (2.5,0.5) {\large $t^2$};
\node[fill opacity = 1] at (0.5,1.5) {\large $q$};
\node[fill opacity = 1] at (1.5,0.5) {\large $t$};
\end{tikzpicture}
\end{center}
\end{ex}

\subsection{Geometric representation}\label{S:ga}
The localized $\T$-equivariant $K$-theory $$K_\T(\PFH_{n,n-k})_\mathrm{loc} = K_\T(\PFH_{n,n-k})\otimes_{\Z[q^{\pm 1},t^{\pm 1}]} \K $$
is a $\K$-vector space with basis given by the classes $[I_{\mu,w}]$ of skyscraper sheaves at the $\T$-fixed points $I_{\mu,w}\in\PFH_{n,n-k}$. The space $\bigoplus_{n\ge 0}\bigoplus_{0\le k\le n}K_\T(\PFH_{n,n-k})_\mathrm{loc}$ affords an action of the algebra $\Atq$ by geometrically-defined operators \cite{CGM}. Here we recall only the action of the $\At$-generators on the fixed-point basis:
\begin{align}
\dd_+ \cdot[I_{\mu,w}] &= -t^k \sum_{\substack{\square\in A(\mu)\\x:=\wt(\square)}} xd_{\mu+x,\mu}\prod_{i=1}^k \frac{x-qw_i}{x-qtw_i}[I_{\mu+x,xw}]\\
\label{dminusformula}\dd_-\cdot [I_{\mu,wx}] &=[I_{\mu,w}]\\
\TT_i\cdot [I_{\mu,w}] &= \frac{(t-1)w_{i+1}}{w_i-w_{i+1}}[I_{\mu,w}]+\frac{w_i-tw_{i+1}}{w_i-w_{i+1}}[I_{\mu,s_i(w)}]
\end{align}
where $A(\mu)$ is the set of $\mu$-addable boxes (relative to $\dg'(\mu)$), $\mu+x=\mu+\square$ for $x=\wt(\square)$ with $\square\in A(\mu)$ is the partition obtained by adding $\square$ to $\mu$, and $d_{\mu+x,\mu}$ is the Macdonald-Pieri coefficient
$$
d_{\mu+x,\mu} = \prod_{\square\in R_{\mu+x,\mu}}\frac{t^{\armm_\mu(\square)}-q^{\legg_\mu(\square)+1}}{t^{\armm_\mu(\square)+1}-q^{\legg_\mu(\square)+1}}\prod_{\square\in C_{\mu+x,\mu}}\frac{t^{\armm_\mu(\square)+1}-q^{\legg_\mu(\square)}}{t^{\armm_\mu(\square)+1}-q^{\legg_\mu(\square)+1}}
$$
where $R_{\mu+x,\mu}$ resp. $C_{\mu+x,\mu}$ is the set of boxes in the row resp. column of $\dg'(\mu)$ corresponding to $\square\in A(\mu)$ such that $x=\wt(\square)$.
We also have
\begin{align}
\label{Tinverse}
\TT_i^{-1}[I_{\mu,w}] &= \frac{(1-t^{-1})w_{i}}{w_i-w_{i+1}}[I_{\mu,w}]+\frac{t^{-1}w_i-w_{i+1}}{w_i-w_{i+1}}[I_{\mu,s_i(w)}].
\end{align}

\subsection{Isomorphism}

As shown in \cite{CGM}, there exists a unique isomorphism of $\At$-modules
$$ \bigoplus_{n\ge 0}\bigoplus_{0\le k\le n} K_\T(\PFH_{n,n-k})_{\mathrm{loc}} \overset{\Phi}{\longrightarrow} \bigoplus_{k\ge 0} \Lambda\otimes\K[y_1,\dotsc,y_k] $$
sending $I_{\varnothing,()}\in K_\T(\PFH_{0,0})_{\mathrm{loc}}$ to $1\in\Lambda=V_0$. Moreover, it is shown that $\Phi([I_\mu])=\tH_\mu$ for any $I_\mu\in \Hilb_n = \PFH_{n,n}$, where $\tH_\mu$ is the modified Macdonald symmetric function. We note that the $\At$-equivariance of $\Phi$ can be upgraded to $\Atq$ by means of an involution $\mathcal{N}$ on $\Atq$ (see \cite{CGM} for details).

\section{Conjecture}

\subsection{Bijections between indexing sets}\label{S:Bij}
We recall a bijection introduced in \cite[Proof of Theorem 7.0.1]{CGM}. Starting with a pair $(\mu,w)$ indexing a $\T$-fixed point in $\PFH_{n,n-k}$, we define $\phi(\mu,w)=(\lambda \, \mid \, \gamma)$ with $\lambda\in\Y$ and $\gamma\in(\Z_{\ge 0})^k$ as follows:
\begin{itemize}
\item Let $\gamma_i=l'(\square_i)$
for $1\le i\le k$, where $\square_i$ is the box of $\mu$ with weight $w_i$; note that $l'(\square_i)+1$ is the height of the column in $\mu$ containing $\square_i$. 

\item Form a partition $\kappa$ by removing all columns from $\mu$ which contain a box $\square_i$ for $1\le i\le k$. Let $\lambda=\kappa^t$.
\end{itemize}

It is easy to see that $\phi$ defines a bijection between the set of $(\mu,w)$ indexing fixed points in $\PFH_{n,n-k}$ and the subset $\lambdagamma$ in  $\Y\times(\Z_{\ge 0})^k$ such that $|\la|+|\ga|=n-k$. The inverse is given is as follows. Starting with $\lambdagamma$, we reconstruct $(\mu,w)$ by the following steps:
\begin{itemize}
\item Let $\kappa=\lambda^t$. 

To form $\mu$, insert $k$ labeled columns into $\kappa$, with heights $\gamma_1+1,\dotsc,\gamma_k+1$ and labels $1,\dotsc,k$, respectively, so that each labeled column lies to the right of any unlabeled column in $\mu$ which has the same height and the labeled columns with the same height have increasing labels from right to left.
 
\item For each $1\le i\le k$, define $w_i$ to be the weight of the box at the top (French) or bottom (English) of the column labeled $i$ in $\mu$. That is
\begin{align*}
w_i &= q^{\gamma_i} t^{c_i}
\end{align*}
where
\begin{align}\label{E:c-def}
c_i &= |\{j : \lambda_j\ge\gamma_i+1\}|+|\{j>i : \gamma_j=\gamma_i\}|+|\{j : \gamma_j>\gamma_i\}|.
\end{align}
\end{itemize}

\subsection{Conjecture}

We are now ready to formulate:

\begin{conj}\label{C:main}
Under the $\Atq$-module isomorphism 
$$ \bigoplus_{n\ge 0}\bigoplus_{0\le k\le n} K_\T(\PFH_{n,n-k})_{\mathrm{loc}} \overset{\Phi}{\longrightarrow} \bigoplus_{k\ge 0}\Lambda\otimes\K[y_1,\dotsc,y_k], $$
the fixed-point classes are sent to modified partially-symmetric Macdonald polynomials, as follows:
\begin{align}\label{E:Conj-Mapping}
(-1)^{|\mu|} q^{n(\mu)}t^{n(\mu')} [I_{\mu,w}] \mapsto t^{n(\mathrm{sort}(\lambda,\gamma))+|(\lambda|\gamma)|}\mathcal{J}_{(\lambda|\gamma)}^{w_0}\Big(\frac{X}{t^{-1}-1}\,\big|\,  y\Big)^*
\end{align}
where $\phi(\mu,w)=\lambdagamma$, $\cJ^{w_0}_{\lambdagamma}(X\midd y)$ is defined in \eqref{E:J-w0}, and $f(X\midd y)^*=f(X\midd y)|_{t\mapsto t^{-1}}$.
\end{conj}

Here $f(X\midd y)\mapsto f\big( \frac{X}{1-t} \midd y\big)$ denotes the $\K[y_1,\dotsc,y_k]$-linear automorphism of $V_k$ given by the plethysm
$$ p_i\Big(\frac{X}{1-t}\Big) = \frac{p_i(X)}{1-t^i}$$
on the first tensor factor $\Lambda$. Conjecture~\ref{C:main} is known to hold on the $k=0$ summands, as shown in \cite{CGM}; the right-hand side of \eqref{E:Conj-Mapping} is the modified Macdonald symmetric function $\tH_\la$ in this case.

Following \cite{CGM}, let us introduce the notation 
\begin{align}\label{E:H-def}
H_{\mu,w}=(-1)^{|\mu|} q^{n(\mu)}t^{n(\mu')} [I_{\mu,w}].
\end{align}
We also define 
$$ \tH_{\lambdagamma} = \tH_{\lambdagamma}(X\midd y) = t^{n(\mathrm{sort}(\lambda,\gamma))+|(\lambda|\gamma)|}\mathcal{J}_{(\lambda|\gamma)}^{w_0}\Big(\frac{X}{t^{-1}-1}\,\big|\,  y\Big)^* $$
and call these the \textit{modified partially-symmetric Macdonald functions}.
The assertion of \eqref{E:Conj-Mapping} is then simply that $\Phi(H_{\mu,w})=\tH_{\lambdagamma}$.

The goal of this paper is to begin a proof of this conjecture by matching matrix elements for the action of $\Atq$ on both sides. More precisely, our main result is as follows:

\begin{thm}\label{T:main}
Let $\Phi'$ be the linear map defined by assignment \eqref{E:Conj-Mapping} on basis elements. Then $\Phi'$ respects the action of the following elements of $\Atq$:
\begin{enumerate}
\item $\TT_1,\TT_2,\dotsc$,
\item for any $k\ge 0$, the element $\dd_-\TT_{k}^{-1}\dotsm \TT_1^{-1} \dd_+\id_k$, which acts on $\Lambda\otimes\K[y_1,\dotsc,y_k]$ as multiplication by the elementary symmetric function $e_1(X)$.
\end{enumerate}
\end{thm}

Since $\Phi$ and $\Phi'$ agree at $k=0$, and in particular on the element $I_{\varnothing,()}$, it suffices to verify that $\Phi'$ is $\Atq$-equivariant in order to establish $\Phi=\Phi'$. In fact, it suffices to establish equivariance for $\At$, the subalgebra generated by $\dd_+,\dd_-$ and $\TT_1,\dotsc$. However, it is difficult to compute the action of $\dd_\pm$ directly in terms our elements on the right-hand side of \eqref{E:Conj-Mapping}. We will complete the proof of Conjecture~\ref{C:main} in forthcoming work with M. Bechtloff Weising using a variant of this idea---namely, by combining our Theorem~\ref{T:main} with his results on limit Cherednik operators \cite{BW}.

An explicit example is given in the next subsection. We also show for this example that $\Phi'$ respects the action of the $y_i$. It should be possible to prove this in general using the same type of argument as the one we give for $e_1(X)$, based on results of \cite{G}, though we will not need it in the sequel. 

Let us comment on the normalization given on the right-hand side of \eqref{E:Conj-Mapping}. The factor $t^{n(\mathrm{sort}(\lambda|\gamma))+|(\la|\ga)|}$ is required for the $e_1(X)$-equivariance, as our proof will show. The factor $t^{-\ell(w_0)}$ in \eqref{E:J-w0} is also forced upon us, as we explain below (see \S\ref{S:Normalization}).

\subsection{Example}

One may compute directly that
\begin{align*}
e_1(X)\cJ_{(\varnothing|0,1)}(X\midd y) &= \frac{1}{1-qt}\cJ_{(2|0,0)}(X\midd y)+\frac{1-q}{(1-t)(1-qt)}\cJ_{(1|0,1)}(X \midd y).
\end{align*}
This implies
\begin{align*}
e_1(X)\cJ_{(\varnothing|0,1)}^{w_0}\Big(\frac{X}{t^{-1}-1}\midd y\Big)^* &= \frac{t-1}{1-qt^{-1}}\cJ_{(2|0,0)}^{w_0}\Big(\frac{X}{t^{-1}-1}\midd y\Big)^*+\frac{1-q}{1-qt^{-1}}t\cJ_{(1|0,1)}^{w_0}\Big(\frac{X}{t^{-1}-1}\midd y\Big)^*
\end{align*}
and hence
\begin{align*}
e_1(X)\tH_{(\varnothing|0,1)} &= \frac{t-1}{t-q}\tH_{(2|0,0)}+\frac{1-q}{t-q}\tH_{(1|0,1)}.
\end{align*}
We note that the operator $w_0T_{w_0}$ from the definition \eqref{E:J-w0} of $\cJ^{w_0}_{(\la|\ga)}$ commutes with multiplication by $e_1(X)$.

For $\mu=(2,1)$ and $w=(t,q)$, which satisfies $\phi(\mu,w)=(\varnothing|0,1)$, we have
\begin{align*}
\dd_+\cdot [I_{(2,1),(t,q)}] &= -t^2 t^2 \frac{(t-q^2)(1-q)}{(t^2-q^2)(t-q)}\frac{(t^2-qt)(t^2-q^2)}{(t^2-qt^2)(t^2-q^2t)}[I_{(3,1),(t^2,t,q)}]\\
&= -t^2 [I_{(3,1),(t^2,t,q)}]\\
\TT_2^{-1}\TT_1^{-1} \dd_+ \cdot [I_{(2,1),(t,q)}] 
&=-t^2 \frac{t-1}{t-q}[I_{(3,1),(t^2,t,q)}]-t^2\frac{1-q}{t-q}[I_{(3,1),(t^2,q,t)}]
\end{align*}
and hence
\begin{align*}
\dd_-\TT_2^{-1}\TT_1^{-1}\dd_+ \cdot H_{(2,1),(t,q)}
&= \frac{t-1}{t-q}H_{(3,1),(t^2,t)}+\frac{1-q}{t-q}H_{(3,1),(t^2,q)}
\end{align*}
where
$H_{(2,1),(t,q)}=(-1)^3 qt [I_{(2,1),(t,q)}],
H_{(3,1),(t^2,t)}=(-1)^4qt^3[I_{(3,1),(t^2,t)}],$ and $H_{(3,1),(t^2,q)}=(-1)^4qt^3[I_{(3,1),(t^2,q)}].$

We note that if we had started with $\cJ_{(\varnothing|1,0)}$, the Pieri formula would have three terms; the index $(\varnothing|1,0)$ corresponds to $\la=(2,1)$, $w=(q,t)$.


Let us also consider the operator $\yy_2$ acting on $H_{(2,1),(q,t)}$. Using $\yy_2=\TT_1^{-1}\varphi$ where $(t-1)\varphi=[\dd_+,\dd_-]$, we find by \cite[(6.2.1)]{CGM} that
\begin{align*}
y_2\cdot H_{(2,1),(q,t)} = \frac{1}{t-q}H_{(2,2),(qt,q)} + \frac{(t-1)t}{(q-t)(t^2-q)}H_{(3,1),(t^2,q)} + \frac{1}{q-t^2}H_{(3,1),(q,t^2)}.
\end{align*}
On the other hand, we compute $T_{w_0}^{-1}w_0 y_2w_0T_{w_0}=T_1^{-1}y_1T_1$ on $\cJ_{(\varnothing|1,0)}$ as
\begin{align*}
T_1^{-1}y_1T_1 \cJ_{(\varnothing|1,0)} = \frac{t^{-1}}{1-qt}\cJ_{(\varnothing|1,1)}-\frac{1-t}{(1-qt)(1-qt^2)}t^{-1}\cJ_{(1|0,1)}-\frac{1}{1-qt^2}\cJ_{(1|1,0)}.
\end{align*}
This implies that
$y_2 \cdot \cJ_{(\varnothing|1,0)}^{w_0}\big(\frac{X}{t^{-1}-1}\midd y\big)^*$
is the sum of the terms
\begin{align*}
\frac{1}{t-q}\cdot t^2\cJ_{(\varnothing|1,1)}^{w_0}\Big(\frac{X}{t^{-1}-1}\midd y\Big)^*\\
\frac{(t-1)t}{(q-t)(t^2-q)}\cdot t^2 \cJ_{(1|0,1)}^{w_0}\Big(\frac{X}{t^{-1}-1}\midd y\Big)^*\\
\frac{1}{q-t^2}\cdot t^2\cJ_{(1|1,0)}^{w_0}\Big(\frac{X}{t^{-1}-1}\midd y\Big)^*
\end{align*}
and therefore that
\begin{align*}
y_2 \cdot \tH_{(\varnothing|1,0)} = \frac{1}{t-q}\tH_{(\varnothing|1,1)} + \frac{(t-1)t}{(q-t)(t^2-q)}\tH_{(1|0,1)} + \frac{1}{q-t^2}\tH_{(1|1,0)},
\end{align*}
in exact agreement with the computation of $y_2\cdot H_{(2,1),(q,t)}$, according to the assignment in our conjecture. (We note that multiplication by $y_2$ on $\cJ_{(\varnothing|1,0)}$ does \textit{not} give matching coefficients; the twist by $w_0T_{w_0}$ is essential.)

\subsection{Normalization}\label{S:Normalization}

The reader may wonder if the operator $t^{-\ell(w_0)}w_0 T_{w_0}$ from \eqref{E:J-w0} is correctly normalized, since one could multiply this by any nonzero scalar (depending on $k$) to achieve the same effect on operators via conjugation. Let us briefly justify this choice, by considering the case of $\mu=(k)$ and $w=(t^{k-1},\dotsm,t,1)$, which satisfies $\phi(\mu,w)=(\varnothing|0,\dotsc,0)$.

We have $\cJ_{(\varnothing|0,\dotsc,0)}=1$ and $$
\tH_{(\varnothing|0,\dotsc,0)}=\cJ^{w_0}_{(\varnothing|0,\dotsc,0)}=1,
$$
thanks in particular to the factor $t^{-\ell(w_0)}$ in \eqref{E:J-w0}. 

A necessary condition for Conjecture~\ref{C:main} to be true is that $$\dd_-\cdot\tH_{(\varnothing|0,\dotsc,0)} = \tH_{(1|0,\dotsc,0)}\in V_{k-1},$$ since $\dd_-\cdot H_{(k),(t^{k-1},\dotsc,1)} = H_{(k),(t^{k-1},\dotsc,t)}$.

We clearly have $\dd_-\cdot \tH_{(\varnothing|0,\dotsc,0)} = \dd_-\cdot 1= e_1(X)$. Also, $\cJ^{w_0}_{(1|0,\dotsc,0)}$ is easily seen to be $(1-t)e_1(X)$. Hence $\tH_{(1|0,\dotsc,0)}=e_1(X)$, and this shows that our normalization is consistent with Conjecture~\ref{C:main}.

\section{Demazure-Lusztig operators}

Our goal in this section is to prove assertion (1) from  Theorem~\ref{T:main}. We assume throughout this section that $(\mu,w)$ indexes a fixed point in $\PFH_{n,n-k}$ and that $\phi(\mu,w)=\lambdagamma$. Also let $m=n-k$ and regard $\la$ as an element of $\Y_m\subset(\Z_{\ge 0})^m$.

\subsection{Special arms and legs}

Suppose $\gamma_i>\gamma_{i+1}$ and let $u=(m+i,\gamma_{i+1}+1)\in\dg(\la^-|\ga)$. Then we have
$\gamma_i-\gamma_{i+1}=l(u)+1$
and, using \eqref{E:c-def},
\begin{align*}
c_{i+1}-c_i &= |\{j : \gamma_i\ge \lambda_j\ge \gamma_{i+1}+1\}|+|\{j<i:\gamma_i\ge \gamma_j\ge\gamma_{i+1}+1\}|\notag\\
&\qquad +|\{j>i+1:\gamma_i>\gamma_j\ge \gamma_{i+1}\}|+1\\
&= a(u)\notag
\end{align*}
with arms and legs taken in $\dg(\la^-|\ga)$.
Therefore, under the assumption $\gamma_i>\gamma_{i+1}$, we have
\begin{align}\label{E:wu}
\frac{w_i}{w_{i+1}} = q^{l(u)+1}t^{-a(u)}.
\end{align}

\subsection{Converting $T_i$ to $\TT_i$}\label{S:T}
Continuing to assume $\gamma_i>\gamma_{i+1}$, we now appeal to \eqref{E:TJ} to see that
\begin{align*}
T_{i}^* \cJ_{(\lambda|\gamma)}(X\midd y)^*
&= \frac{q^{-(l(u)+1)}t^{a(u)+1}-1}{q^{-(l(u)+1)}t^{a(u)}-1}\cJ_{(\lambda|s_i(\gamma))}(X\midd y)^*+\frac{t-1}{q^{l(u)+1}t^{-a(u)}-1}\cJ_{(\lambda|\gamma)}(X\midd y)^*
\end{align*}
where
\begin{align*}
T_{i}^* &= s_i + \frac{(1-t)y_{i+1}}{y_{i+1}-y_i}(1-s_i) = \frac{y_i-ty_{i+1}}{y_i-y_{i+1}}s_i + \frac{(1-t)y_{i+1}}{y_{i+1}-y_{i}}
\end{align*}
for $1\le i \le k-1$. More generally, we let $T_w^* = t^{\ell(w)} (* \circ \, T_w \circ *)$ where $*$ is the automorphism sending $t\mapsto t^{-1}$ and fixing everything else. Then $$ w_0 T_{w_0}^* T_i^* (w_0 T_{w_0}^*)^{-1} $$ agrees with the operator $\TT_i$ from the polynomial representation of $\Atq$; this follows from the relation $T_{w_0} T_i T_{w_0}^{-1} = T_{k-i}$ (see, e.g., \cite[(3.1.8)]{M}). Hence we have
\begin{align*}
\TT_i \cdot \tH_{(\lambda|\gamma)}
&= \frac{q^{-(l(u)+1)}t^{a(u)+1}-1}{q^{-(l(u)+1)}t^{a(u)}-1}\tH_{(\lambda|s_i(\gamma))}+\frac{t-1}{q^{l(u)+1}t^{-a(u)}-1}\tH_{(\lambda|\gamma)}.
\end{align*}

On the other hand, by \eqref{E:wu}, we see that
\begin{align*}
\TT_i\cdot [I_{\mu,w}] &= \frac{(t-1)w_{i+1}}{w_i-w_{i+1}}[I_{\mu,w}]+\frac{w_i-tw_{i+1}}{w_i-w_{i+1}}[I_{\mu,s_i(w)}]\\
&= \frac{t-1}{q^{l(u)+1}t^{-a(u)}-1}[I_{\mu,w}]+\frac{q^{-(l(u)+1)}t^{a(u)+1}-1}{q^{-(l(u)+1)}t^{a(u)}-1}[I_{\mu,s_i(w)}].
\end{align*}
This completes the proof of Theorem~\ref{T:main}(1) in this case when $\ga_i>\ga_{i+1}$. Using the Hecke algebra action, the case $\ga_i<\ga_{i+1}$ automatically follows. When $\ga_i=\ga_{i+1}$, the assertion follows immediately from $T_i \cJ_{\lambdagamma} = t \cJ_{\lambdagamma}$  (see, e.g., \cite[(3.3.40)]{Ch}) and $w_i=tw_{i+1}$, so we are done.

\section{Pieri formulas}\label{ElementaryMultiplication}

In this section, we prove assertion (2) from Theorem~\ref{T:main}, which is significantly more difficult than assertion (1). We use a Pieri formula for the partially-symmetric integral forms $\cJ_{\lambdagamma}$ established in \cite{G} using interpolation polynomials.

We make the following notational adjustments in order to match \cite{G} as closely as possible: we work in $\PFH_{N,N-k}$ and use $(\xi,w)$ instead of $(\mu,w)$ to index a typical fixed point in this space, so $|\xi|=N$ and the length of $w$ is $k$. Also, when $\phi(\xi,w)=\lambdagamma$, we may regard $\la$ as an element of $\Y_m\subset (\Z_{\ge 0})^m$  and $\lambdagamma$ as an element of $(\Z_{\ge 0})^n$, where $m=n-k$ and $n>N$ is arbitrary. Our computations will show that the Pieri formula for $\cJ_{\lambdagamma}$ in $n$ total variables does not depend on $n$ in this range.

To state formulas for $e_1(X) \JJ$, we need several definitions. Starting with the eigenvalues\footnote{These are the eigenvalues of Cherednik-Dunkl operators on the nonsymmetric Macdonald polynomials. Even though we will not use these operators directly in this paper, we still refer to the $\overline{\nu}_i$ as eigenvalues.}, we have
\begin{align}\label{eigenvaluesone}\overline{\nu}_i &= q^{\nu_i} t^{-l'_\nu(i)}, \hspace*{.5cm} 1\leq i \leq n,\\
\label{eigenvaluestwo}l'_\nu(i) &= \# \{j<i \mid \nu_j > \nu_i\} + \# \{j>i \mid \nu_j \geq \nu_i\}.\end{align}
 
 We now take $\lambdagamma$ and build several new compositions. Let $I_1:=\{t_1, \dotsc, t_r\}$ with $ 0 = t_0 < t_1 < \, \dotsm \, < t_r < t_{r+1} = k+1$. Define $\eta$ by choosing some entry of $\lambda$ to be $\lambdatilde_{n-k}$ and then setting,
\begin{align}
    &\eta_i = \gamma_i \hspace*{.55cm}\text{ if } i\notin I_1\nonumber \\
    &\eta_{t_j}=\gamma_{t_{j-1}} \text{ for } 1\leq i \leq r \text{ where } \gamma_{t_0} = \lambdatilde_{n-k}.\label{eta}
\end{align}

Next we build $\mutilde$ by rearranging all entries of $\lambda$ except $\lambdatilde_{n-k}$, which is replaced with a new column of height $\gamma_{t_r}+1$. Choose $\mutilde$ such that all columns of height $\gamma_{t_r}+1$ are to the right, starting with the entry in position $h$, so $\mutilde_{h} = \dotsm = \mutilde_{n-k} = \gamma_{t_r}+1$ , and the other columns are weakly increasing, $\mutilde_1 \leq \, \dotsm \, \leq \mutilde_{h-1}$. Finally, we build the rest of $\lambdatilde$ by setting $(\lambdatilde_1, \dotsc, \lambdatilde_{n-k-1}) := (\mutilde_1, \dotsc, \mutilde_{n-k-1})$.

We also require that $I_1$ satisfy the following two properties:
\begin{itemize}
    \item $\eta_j \neq \eta_{t_u}$ for any $u \in \{1,\dotsc, r\}$ and $j \in \{t_{u-1}+1,\dotsc, t_{u}-1\}$, and
    \item $\eta_j \neq \mutilde_{n-k}-1$ for any $j \in \{t_r+1, \dotsc, k\}$.
\end{itemize}

Such an $I_1$ is called maximal with respect to $\mueta$. This can be used to define the support of the expansion, \[ \MMM:= \{ \mueta \, \mid \, \mu \in S_m(\mutilde), \, \mu_1 \geq \dotsm \geq \mu_m, \, \eta \text{ is as defined above for some maximal } I_1 \} \]

\begin{thm}[\cite{G}, Theorem 5.7]\label{Jpieri}
One has an expansion
\[ e_1(X) \JJ = \sum_{\mueta \in \MMM}  \AAA \mathcal{J}_{\mueta} \]
with coefficients given as follows:
\begin{equation}\label{simplifiedpieri} \AAA = \prod_{ \square \in \dg_{(\lambdatilde_{n-k}+1)}(\lambda^-)}\frac{t-q^{\ell_{\lambdamingamma}(\square)+1}t^{a_{\lambdamingamma}(\square)+1}}{1-q^{\ell_{\lambdamingamma}(\square)+1}t^{a_{\lambdamingamma}(\square)+1}} \cdot \jtwo \cdot p'_{I_1} \cdot \left(\frac{1}{1-t}\right) \cdot (q^{-\mutilde_h+1}) \cdot \overline{\lambdatildegamma}_{n-k}\end{equation}
where 
\[ p'_{I_1} := \left(\frac{(t-1)q^{-1} \overline{\mutilde}_{\hIone}}{q^{-1} \overline{\mutilde}_{\hIone}-\overline{\eta}_{t_r}}\right)  \prod_{u=1}^{r-1} \left(\frac{(t-1) \overline{\eta}_{t_{u+1}}}{ \overline{\eta}_{t_{u+1}}-\overline{\eta}_{t_{u}}} \right)  \prod_{j=t_r+1}^{k} \left( \frac{(q^{-1}t) \overline{\mutilde}_{\hIone} -  \overline{\eta}_j}{q^{-1} \overline{\mutilde}_{\hIone} -  \overline{\eta}_j}\right)\prod_{u=1}^{r} \prod_{j=t_{u-1}+1}^{t_{u}-1} \left( \frac{t\overline{\eta}_{t_u} -  \overline{\eta}_j}{ \overline{\eta}_{t_{u}} -  \overline{\eta}_j} \right),\]
for $0 = t_0 < t_1 < \dotsm < t_r < t_{r+1}  = k+1$, and
\[\jtwo  = \prod\limits_{\square \in \mathrm{C}_{\lambda^-}}\frac{1-q^{\leg{\lambdamingamma}}t^{\armone{\lambdamingamma}+1} }{ 1-q^{\leg{\lambdamingamma}+1}t^{\armone{\lambdamingamma}+1+m_{j-1}(\eta)}},\]
where $\mathrm{C}_{\lambda^-}$ is the rightmost column in $\dg(\lambda^-)$ of height $\lambdatilde_{n-k}$ and $m_{j-1}(\eta)$ is the number of columns in $\eta$ of height $j-1$.

\end{thm}

\subsection{$p'_{I_1}$ vs. coefficient of $\TT_{k}^{-1} \dotsm \TT_{1}^{-1}$}

To compare $p'_{I_1}$ and the coefficient resulting from $\TT_{k}^{-1} \dotsm \TT_{1}^{-1}$, we first rewrite $p'_{I_1}$ using the eigenvalues $\overline{\gamma}_i$ replacing their respective $\overline{\eta}_i$:
\begin{equation}
    p'_{I_1} =  \frac{(t-1)\overline{\gamma}_{t_1}}{\overline{\gamma}_{t_1} - \overline{\eta}_{t_{1}}} \cdot \prod_{u=1}^{r-1} \frac{(t-1)\overline{\gamma}_{t_{u+1}}}{\overline{\gamma}_{t_{u+1}} - \overline{\gamma}_{t_{u}}}  \cdot \prod_{j=1}^{t_1-1}\frac{t \overline{\eta}_{t_{1}}-\overline{\gamma}_j}{\overline{\eta}_{t_1}-\overline{\gamma}_j} \cdot \prod_{u=1}^{r} \prod_{j=t_{u}+1}^{t_{u+1}-1} \frac{t \overline{\gamma}_{t_{u}}-\overline{\gamma}_j}{\overline{\gamma}_{t_u}-\overline{\gamma}_j}
\end{equation}
Comparing the formulas for the eigenvalues in \eqref{eigenvaluesone} and \eqref{eigenvaluestwo} to the weights in \eqref{E:c-def} gives a simple conversion,
\[(\overline{\gamma}_i)_{t\mapsto t^{-1}} = w_i \qquad \text{ and } \qquad (\overline{\eta}_{t_1})_{t\mapsto t^{-1}} = x.\]
This allows us to rewrite,
\begin{equation}\label{HHLpIone}
    \left[p'_{I_1}\right]^* = \frac{(t^{-1}-1) w_{t_1}}{w_{t_1}-x} \cdot \prod_{u=1}^{r-1} \frac{(t^{-1}-1)w_{t_{u+1}}}{w_{t_{u+1}}-w_{t_u}} \cdot \prod_{j=1}^{t_1-1} \frac{t^{-1}x-w_j}{x-w_j} \cdot \prod_{u=1}^r \prod_{j=t_u+1}^{t_{u+1}-1} \frac{t^{-1} w_{t_u} - w_j}{w_{t_u}-w_j},
\end{equation}
where $*$ is the automorphism of $\K$ sending $q\mapsto q$ and $t\mapsto t^{-1}$.

We want to match the coefficients in the expansion of $e_1(X) $  with $\dd_-\TT_k^{-1}\dotsm \TT_1^{-1} \dd_+$. Since $\dd_+ H_{\xi,w}$ is a linear combination of some $H_{\xi+x,xw}$, the product $\TT_k^{-1}\dotsm \TT_1^{-1}$ acts on $H_{\xi+x,xw}$ where $xw=(x,w_1,\dots,w_k)$.

We will use,
\begin{equation}\label{LeftRightTi}
\TT_i^{-1}(H_{\xi,w}) = \frac{(1-t^{-1})w_{i}}{w_i-w_{i+1}}H_{\xi,w}+\frac{t^{-1}w_i-w_{i+1}}{w_i-w_{i+1}}H_{\xi,s_i(w)}.
\end{equation}

\begin{lem} Let $v = \reversesproduct$. The full action $\TT_k^{-1} \dotsm \TT_1^{-1} H_{\xi+x,xw}$ can be written
\begin{equation}\label{CGMTProduct}
    \TT_k^{-1} \dotsm \TT_1^{-1} H_{\xi+x,xw} = \sum_{I_1\subseteq [1,k]} \tildepIone \cdot H_{\xi+x,v(xw)}
\end{equation}
where
\begin{equation}\label{CGMpIone} \tildepIone =  \frac{(1-t^{-1})x}{x-w_{t_1}} \prod_{u=1}^{r-1} \frac{(1-t^{-1})w_{t_u}}{w_{t_u}-w_{t_{u+1}}} \cdot \prod_{j=1}^{t_1-1} \frac{t^{-1}x-w_j}{x-w_j} \cdot \prod_{u=1}^{r} \prod_{j=t_{u}+1}^{t_{u+1}-1} \frac{t^{-1}w_{t_{u}} - w_j}{w_{t_u} - w_j}. \end{equation}
\end{lem}

\begin{proof}
     In the expansion of $\TT_k^{-1} \dotsm \TT_1^{-1} H_{\xi+x,xw}$, an additive term can be identified by a choice for each $i\in [1,k]$ of either the left rational function in \eqref{LeftRightTi}, which does not include the transposition $s_i$, or the right function which does. Fix a particular summand in the expansion, and denote by $v:=\reversesproduct$  the associated product of transpositions. The hats indicate absence of that simple transposition. We choose the left function of \eqref{LeftRightTi} for $\ell \in I_1 := \{t_1, \dotsc, t_k\}$, and the right function for $\ell \notin I_1$. Call $f_\ell$ the coefficient contributed by each $\TT^{-1}_\ell$ acting on $\TT^{-1}_{\ell-1} \dotsm \TT^{-1}_1 H_{\xi+x,xw}$, so in total,
     \[\TT_k^{-1} \dotsm \TT_1^{-1} H_{\xi+x,xw} = (f_k \dotsm f_1) \reversesproduct H_{\xi+x,xw} = f_k \dotsm f_1 H_{\xi+x,v(xw)} \]
    Then the following are computations of the coefficients $f_\ell$.
    
    \underline{Case 1:} Let $\ell < t_1$, and define $\sigma := s_{\ell-1} \dotsm s_1$. 
     \[ \TT_\ell^{-1} \left( \TT^{-1}_{\ell-1} \dotsm \TT^{-1}_1 H_{\xi+x,xw} \right)  =  \frac{t^{-1} (\sigma(xw))_\ell - (\sigma(xw))_{\ell+1}}{(\sigma(xw))_\ell - (\sigma(xw))_{\ell+1}} s_\ell  \left( f_{\ell-1} \dotsm f_1 H_{\xi+x,\sigma(xw)}\right).\]
    The permutation $\sigma$ moves $x$ into position $\ell$, so $((\sigma(xw))_\ell=x$, and $(\sigma(xw))_{\ell+1} = w_\ell$. So explicitly,
    \[f_\ell = \frac{t^{-1} x - w_{\ell}}{x-w_\ell}.\]
    
    \underline{Case 2:}  Let $\ell=t_1$ and $\sigma := s_{t_1-1} \dotsm s_1$. Then like in Case 1,
    \begin{align*} \TT^{-1}_{t_1}\left( \TT^{-1}_{t_1-1} \dotsm \TT^{-1}_1 H_{\xi+x,xw} \right)  &= \frac{(1-t^{-1})\cdot (\sigma(xw))_{t_1}}{(\sigma(xw))_{t_1} - (\sigma(xw))_{t_1+1}} \left( f_{t_1-1} \dotsm f_1 H_{\xi+x,\sigma(xw)}\right)\\
    f_{t_1} &= \frac{(1-t^{-1}) x}{x-w_{t_1}}\end{align*}
    
    \underline{Case 3:} Let $\ell=t_i$ with $i > 1$, and write $\sigma := (s_{t_i-1} \dotsm s_{t_{i-1}+1}) \widehat{s}_{t_{i-1}} \dotsm s_1$.
    \[ \TT^{-1}_{t_i} \left( \TT^{-1}_{t_i-1} \dotsm \TT^{-1}_1 H_{\xi+x,xw} \right) = \frac{(1-t^{-1})\cdot(\sigma(xw))_{t_i}}{(\sigma(xw))_{t_i} - (\sigma(xw))_{t_i+1}} \left( f_{t_i-1} \dotsm f_1 H_{\xi+x,\sigma(xw)}\right).\]
    The missing transposition $\widehat{s}_{t_{i-1}}$ in $\sigma$ leaves $w_{t_{i-1}}$ in the $(t_{i-1}+1)$st position, and the product $s_{t_i-1}\dotsm s_{t_{i-1}+1}$ moves $w_{t_{i-1}}$ into the $(t_i)$th position. So $(\sigma(xw))_{t_i} = w_{t_{i-1}}$. And the $(t_i+1)$st position is untouched, so $(\sigma(xw))_{\ell+1}=w_{t_i}$. As a result,
    \[ f_{t_i} = \frac{(1-t^{-1}) w_{t_{i-1}}}{w_{t_{i-1}} - w_{t_i}}\]
    
    \underline{Case 4:} Let $t_{i-1} < \ell < t_{i}$ for $i>1$, and again let $\sigma := (s_{\ell-1} \dotsm s_{t_{i-1}+1}) \widehat{s}_{t_{i-1}} \dotsm s_1$. Like in Case 1, we have,
    \[ \TT_\ell^{-1} \left( \TT^{-1}_{\ell-1} \dotsm \TT^{-1}_1 H_{\xi+x,xw} \right)  = \frac{t^{-1} (\sigma(xw))_\ell - (\sigma(xw))_{\ell+1}}{(\sigma(xw))_\ell - (\sigma(xw))_{\ell+1}} s_\ell \left( f_{\ell-1} \dotsm f_1 H_{\xi+x,\sigma(xw)}\right).\]
    The relevant portion of $\sigma$ is the product $s_{\ell-1} \dotsm s_{t_{i-1}+1}$, which moves $w_{t_i}$ to the $(\ell)$th position, making $(\sigma(xw))_{\ell}=w_{t_{i-1}}$. And the simpler part is $(\sigma(xw))_{\ell+1} = w_{\ell}$, as it has not yet permuted in the $\TT^{-1}_\ell$ coefficient formula. So we obtain
    \[ f_\ell =  \frac{t^{-1}w_{t_{i-1}} - w_\ell}{w_{t_{i-1}} - w_\ell} s_\ell. \]

    Combining all the coefficients, $\tildepIone = f_k \dotsm f_1$, gives the desired formula \eqref{CGMpIone}.
\end{proof}

\begin{rem}
     A priori, the formula \eqref{CGMTProduct} includes all subsets $I_1 \subseteq[1,k]$, which would include those which swap labels of boxes in the same row. That cannot be allowed, as $v(xw)$ would no longer correspond to a flag of ideals. However, any such term will vanish, as adjacent boxes being swapped causes either $(t^{-1}x-w_j)$ or $(t^{-1}w_{t_u}-w_j)$ to be 0. This coincides with the maximality condition on $I_1$ in the support $\MMM$ of $e_1(X) \JJ$.
\end{rem}

\begin{cor}
\begin{equation}\label{pIoneComparison} \tildepIone = \left[p'_{I_1}\right]^* \cdot x \cdot (w_{t_r})^{-1}\end{equation}
\end{cor}
 \begin{proof}
 This comes immediately from comparing terms in \eqref{HHLpIone} with \eqref{CGMpIone}, which only differ in the numerators in the leftmost two rational functions.
 \end{proof}

\subsection{Remaining pieces} First we must address some internal cancellation from terms in ~\cite{CGM}. We recall the formula for $\dd_+ H_{\xi,w}$ and Pieri coefficients $d_{\xi+x,\xi}$ from ~\cite{CGM} here:
 \begin{align*} \dd_+ H_{\xi,w} &= t^k \sum_x d_{\xi+x,\xi} \prod_{i=1}^k \frac{x-qw_i}{x-qtw_i} H_{\xi+x,xw}\\
 d_{\xi+x,\xi} &= \prod_{\square \in R_{\xi+x,\xi}} \frac{t^{\armm_\xi(\square)}-q^{\legg_{\xi}(\square)+1}}{t^{\armm_{\xi}(\square)+1}-q^{\legg_{\xi}(\square)+1}} \prod_{\square \in C_{\xi+x,\xi} } \frac{t^{\armm_\xi(\square)+1}-q^{\legg_\xi(\square)}}{t^{\armm_{\xi}(\square)+1}-q^{\legg_{\xi}(\square)+1}}  \end{align*}
 We split $R_{\xi+x,\xi}$ into two groups,
\begin{align}
    \RS &:= \{\square \in R_{\xi+x,\xi} \, \mid \, \square \text{ is not in a labeled column} \}\\
    \RNS&:=   \{\square \in R_{\xi+x,\xi} \, \mid \, \square \text{ is in a labeled column} \},
\end{align}
which are the boxes in the same row as $x$ that are in `symmetric' and `nonsymmetric' columns respectively. We break up the column labels $[1,k]$ as well,
\begin{align}
    \LL &:= \{ i \in [1,k] \, \mid \, \text{ The column labeled } i \text{ is to the left of } x \}\\
    \RL &:= \{ i \in [1,k] \, \mid \,  \text{ The column labeled } i \text{ is to the right of } x \},
\end{align}
the `left-labeled' and `right-labeled' columns respectively. 

\begin{ex}
The following is the diagram for $H_{\xi,w}=H_{(10,9,6,4,3,2),(q^4t^2,qt^8,q^5t,q^2t^5,t^9)}, $ and $H_{\xi+x,xw} = H_{(10,9,7,4,3,2),(q^2t^6,q^4t^2,qt^8,q^5t,q^2t^5,t^9)}$. Yellow boxes belong to $\RS$, and blue belong to $\RNS$. The labeled columns are $\LL = \{1,3,4\}$ and $\RL = \{2,5\}$.
\end{ex}

\begin{center}
    \begin{tikzpicture}[scale=0.6,fill opacity=.3,every node/.style={scale=0.6}]

\fill[blue,opacity=0.3] (1,2) rectangle (2,3);
\fill[blue,opacity=0.3] (2,2) rectangle (3,3);
\fill[blue,opacity=0.3] (5,2) rectangle (6,3);

\fill[yellow,opacity=0.5] (0,2) rectangle (1,3);
\fill[yellow,opacity=0.5] (3,2) rectangle (4,3);
\fill[yellow,opacity=0.5] (4,2) rectangle (5,3);

\draw (0,0) -- (10,0);
\draw (0,1) -- (10,1);
\draw (0,2) -- (9,2);
\draw (0,3) -- (6,3);
\draw[dashed] (6,3) -- (7,3);
\draw (0,4) -- (4,4);
\draw (0,5) -- (3,5);
\draw (0,6) -- (2,6);
\draw (0,0) -- (0,6);
\draw (1,0) -- (1,6);
\draw (2,0) -- (2,6);
\draw (3,0) -- (3,5);
\draw (4,0) -- (4,4);
\draw (5,0) -- (5,3);
\draw (6,0) -- (6,3);
\draw (7,0) -- (7,2);
\draw[dashed] (7,2) -- (7,3);
\draw (8,0) -- (8,2);
\draw (9,0) -- (9,2);
\draw (10,0) -- (10,1);

\draw (1.5,5.5) circle (13pt);
\node[fill opacity = 1] at (1.5,5.5)  {\small $3$};
\draw (2.5,4.5) circle (13pt);
\node[fill opacity = 1] at (2.5,4.5)  {\small $1$};
\draw (5.5,2.5) circle (13pt);
\node[fill opacity = 1] at (5.5,2.5)   {\small $4$};
\draw (8.5,1.5) circle (13pt);
\node[fill opacity = 1] at (8.5,1.5)  {\small $2$};
\draw (9.5,0.5) circle (13pt);
\node[fill opacity = 1] at (9.5,0.5)  {\small $5$};

\node[fill opacity = 1] at (6.5,2.5) {\large $x$};

\end{tikzpicture}
\end{center}

\begin{lem}
If $i\in \LL$ and $\square$ is the unique box in $R_{\xi+x,\xi}$ in the column labeled $i$, then
\[ \left(\frac{t^{\armm_\xi(\square)} - q^{\legg_\xi(\square)+1}}{t^{\armm_{\xi}(\square)+1} - q^{\legg_{\xi}(\square)+1}} \right)\cdot \left(\frac{x-qw_i}{x-qtw_i}\right) = t^{-1}.\]
\end{lem}

\begin{proof}
We will write $u$ to indicate the coordinate of $\square$, $u = q^r t^c$. Then we have
    \begin{align}
    \frac{x}{u}=t^{\armm_{\xi+x}(\square)},\qquad 
    \frac{w_i}{u} = q^{\legg_{\xi+x}(\square)}
    \end{align}
    and hence
    \begin{align}
    \frac{x-qw_i}{x-qtw_i} &= \frac{\frac{x}{u}-q\frac{w_i}{u}}{\frac{x}{u}-qt\frac{w_i}{u}} = \frac{t^{\armm_{\xi+x}(\square)}-q^{\legg_{\xi+x}(\square)+1}}{t^{\armm_{\xi+x}(\square)}-tq^{\legg_{\xi+x}(\square)+1}} = t^{-1} \frac{t^{\armm_\xi(\square)+1}-q^{\legg_\xi(\square)+1}}{t^{\armm_\xi(\square)}-q^{\legg_\xi(\square)+1}}
    \end{align}
    since $\armm_{\xi+x}(\square)=\armm_\xi(\square)+1$ and $\legg_{\xi+x}(\square)=\legg_\xi(\square)$.
\end{proof}

Cancelling terms from every nonsymmetric column, we have the resulting equation,
\begin{equation}\label{rowsimplification}
    \prod_{\square \in R_{\xi+x,\xi}} \frac{t^{\armm_\xi(\square)}-q^{\legg_\xi(\square)+1}}{t^{\armm_\xi(\square)+1} - q^{\legg_{\xi}(\square)+1}}  \prod_{i \in \LL} \frac{x-qw_i}{x-qtw_i} = t^{-|\LL|} \prod_{\square \in \RS} \frac{t^{\armm_\xi(\square)}-q^{\legg_\xi(\square)+1}}{t^{\armm_\xi(\square)+1} - q^{\legg_{\xi}(\square)+1}}.
\end{equation}

We now do something similar with the column terms, which behaves differently since a box in the same column as $x$ may have multiple labeled boxes to its right. 

\begin{lem}
Let $\square = (i,j) \in C_{\xi+x,\xi}$, and let $\ell_1 <  \dotsm <  \ell_b \in \RL$ be the labels of all labeled boxes in row $j$, the same row as $\square$. Then we have
\[
    \frac{t^{\armm_\xi(\square)+1}-q^{\legg_\xi(\square)}}{t^{\armm_{\xi}(\square)+1}-q^{\legg_{\xi}(\square)+1}}  \cdot \prod_{i=1}^b \frac{x-qw_{\ell_i}}{x-qtw_{\ell_i}} =  \frac{t^{\armm_{\xi+x}(\square)+1-b} - q^{\legg_{\xi+x}(\square)-1}}{t^{\armm_{\xi+x}(\square)+1} - q^{\legg_{\xi+x}(\square)}}.\]
\end{lem}

\begin{proof}
Since the labels $\ell_1, \dotsc,  \ell_b$ must appear in adjacent columns, with $\ell_1$ the farthest to the right, we can write $w_{\ell_1} = t w_{\ell_2} = \dots = t^{b-1} w_{\ell_b}.$ Therefore we can simplify,
\[\prod_{i=1}^b \frac{x-qw_{\ell_i}}{x-qtw_{\ell_i}}  = \prod_{i=1}^b \frac{x-t^{-i+1}qw_{\ell_1}}{x-t^{-i+2} qw_{\ell_1}} = \frac{x-t^{-b+1}qw_{\ell_1}}{x-t qw_{\ell_1}}\]
Once again writing $u$ for the coordinate of $\square$, we compare arm and leg counts,
\[ \frac{x}{u} = q^{\legg_{\xi+x}(\square)}, \qquad \frac{w_{\ell_1}}{u} = t^{\armm_{\xi+x}(\square)},\]
which allows us to simplify,
\[ \prod_{i=1}^{b}\frac{x-qw_{\ell_i}}{x-qtw_{\ell_i}} = \frac{\frac{x}{u} - t^{-b+1}q \frac{w_{\ell_1}}{u}}{\frac{x}{u} - tq \frac{w_{\ell_1}}{u}} = \frac{q^{\legg_{\xi+x}(\square)} - q t^{\armm_{\xi+x}(\square)+1-b}}{q^{\legg_{\xi+x}(\square)} - q t^{\armm_{\xi+x}(\square)+1}}=\frac{t^{\armm_{\xi+x}(\square)+1-b} - q^{\legg_{\xi+x}(\square)-1}}{t^{\armm_{\xi}(\square)+1} - q^{\legg_{\xi}(\square)}},\]
where in the last equality we use,
\[ \armm_{\xi+x}(\square) = \armm_\xi(\square), \qquad \legg_{\xi+x}(\square) = \legg_{\xi}(\square)+1. \qedhere\]
\end{proof}

Taking all of these terms together, we have the following identity:
\begin{equation}\label{columnsimplification}
    \prod_{\square \in C_{\xi+x,\xi}} \frac{t^{\armm_\xi(\square)+1}-q^{\legg_\xi(\square)}}{t^{\armm_{\xi}(\square)+1}-q^{\legg_{\xi}(\square)+1}}   \cdot \prod_{i\in \RL} \frac{x-qw_i}{x-qtw_i} = \prod_{\square \in C_{\xi+x,\xi}} \frac{t^{\armm_{\xi+x}(\square)+1-\cbox} - q^{\legg_{\xi+x}(\square)-1}}{t^{\armm_{\xi+x}(\square)+1} - q^{\legg_{\xi+x}(\square)}},
\end{equation}
where $\cbox$ is the number of labeled boxes in the same row as $\square$.

\begin{thm}\label{T:Pieri-match}
Consider the expansion
\begin{equation}
    \left(\dd_- \TT_k^{-1} \dotsm \TT_1^{-1} \dd_+\right) H_{\xi,w} = \sum_{x,v} C^{\xi,w}_{\xi+x,(v(xw))'} H_{\xi+x,(v(xw))'},
\end{equation}
where $v$ is the permutation $v=\reversesproduct$, and $(v(xw))'$ indicates the last entry of $v(wx)$ is removed. The coefficients satisfy
\[ C^{\xi,w}_{\xi+x,(v(xw))'} = t^{c_r} \cdot \left[\AAA \cdot (1-t)\right]^*,\]
where $\xi,w$ corresponds to $\lambdagamma$ and $\xi+x,(v(xw))'$ corresponds to $\mueta$ as in Section $\ref{S:Bij}$, the column index of $w_{t_r}$ is $c_r$, and $*$ is the $\mathbb{K}$-automorphism fixing $q$ and sending $t\mapsto t^{-1}$.
\end{thm}

\begin{proof}
Applying \eqref{rowsimplification} and \eqref{columnsimplification}, we have a simplified formula for $\displaystyle{ \dd_+ H_{\xi,w} },$  \[\sum_x  t^{|\RL|} \prod_{\square \in \RS} \frac{t^{\armm_\xi(\square)}-q^{\legg_\xi(\square)+1}}{t^{\armm_\xi(\square)+1} - q^{\legg_{\xi}(\square)+1}} \prod_{\square \in C_{\xi+x,\xi}} \frac{t^{\armm_{\xi+x}(\square)+1-\cbox} - q^{\legg_{\xi+x}(\square)-1}}{t^{\armm_{\xi+x}(\square)+1} - q^{\legg_{\xi+x}(\square)}} H_{\xi+x,xw}.\]

Combining this form with \eqref{CGMpIone} and \eqref{dminusformula} gives a formula for $C^{\xi,w}_{\xi+x,(v(xw))'
}$,

\[  \tildepIone \cdot  t^{|\RL|} \prod_{\square \in \RS} \frac{t^{\armm_\xi(\square)}-q^{\legg_\xi(\square)+1}}{t^{\armm_\xi(\square)+1} - q^{\legg_{\xi}(\square)+1}} \prod_{\square \in C_{\xi+x,\xi}} \frac{t^{\armm_{\xi+x}(\square)+1-\cbox} - q^{\legg_{\xi+x}(\square)-1}}{t^{\armm_{\xi+x}(\square)+1} - q^{\legg_{\xi+x}(\square)}}  \]

On the other side, from \eqref{simplifiedpieri},
\[\AAA \cdot (1-t) = \displaystyle{\prod_{ \square \in \dg_{(\lambdatilde_{n-k}+1)}(\lambda^-)}\frac{t-q^{\ell_{\lambdamingamma}(\square)+1}t^{a_{\lambdamingamma}(\square)+1}}{1-q^{\ell_{\lambdamingamma}(\square)+1}t^{a_{\lambdamingamma}(\square)+1}}} \cdot \jtwo \cdot p'_{I_1} \cdot (q^{-\mutilde_h+1}) \cdot \overline{\lambdatildegamma}_{n-k}.\]
It therefore remains to prove the following correspondences, where $*$ is again the $\mathbb{K}$-\\automorphism fixing $q$ and sending $t\mapsto t^{-1}$.

    \begin{enumerate}
         \item $\displaystyle{\left[\prod_{ \square \in \dg_{(\lambdatilde_{n-k}+1)}(\lambda^-)}\frac{t-q^{\ell_{\lambdamingamma}(\square)+1}t^{a_{\lambdamingamma}(\square)+1}}{1-q^{\ell_{\lambdamingamma}(\square)+1}t^{a_{\lambdamingamma}(\square)+1}}\right]^* = \prod_{\square \in \RS} \frac{t^{\armm_\xi(\square)}-q^{\legg_\xi(\square)+1}}{t^{\armm_\xi(\square)+1} - q^{\legg_\xi(\square)+1}}}$ \\[3mm]

        \item $\displaystyle{\left[\jtwo\right]^* = t^{|\RL|}} \cdot \prod_{\square \in C_{\xi+x,\xi}} \frac{t^{\armm_{\xi+x}(\square)+1-\cbox} - q^{\legg_{\xi+x}(\square)-1}}{t^{\armm_{\xi+x}(\square)+1} - q^{\legg_{\xi+x}(\square)}} $
       
         \item $\displaystyle{\left[p'_{I_1}\cdot  (q^{-\mutilde_h+1})  \cdot \overline{\lambdatildegamma}_{n-k}   \right]^* = \tildepIone \cdot   t^{c_r}}$ \\[3mm]
        \end{enumerate}

Starting with (1), begin with the following:
\begin{align*}\left[\prod_{ \square \in \dg_{(\lambdatilde_{n-k}+1)}(\lambda^-)}\frac{t-q^{\ell_{\lambdamingamma}(\square)+1}t^{a_{\lambdamingamma}(\square)+1}}{1-q^{\ell_{\lambdamingamma}(\square)+1}t^{a_{\lambdamingamma}(\square)+1}}\right]^*   &=    \prod_{ \square \in \dg_{(\lambdatilde_{n-k}+1)}(\lambda^-)}\frac{t^{a_{\lambdamingamma}(\square)}-q^{\ell_{\lambdamingamma}(\square)+1}}{t^{a_{\lambdamingamma}(\square)+1}-q^{\ell_{\lambdamingamma}(\square)+1}}. \end{align*}

For a box $\square = (i,\lambdatilde_{n-k}+1)$ in $ d_{(\lambdatilde_{n-k}+1)}(\lambda)$, notice that $\lambdatilde_{n-k}+1$ is the height of the added box $x$. So this product is over all boxes in $\lambdatilde$, the symmetric part of the diagram, of the same height as $x$, each of which corresponds to a box $\square' \in \RS$. The arm and leg counts for the corresponding boxes are easily compared using the respective definitions to get,
\[ \legg_\xi(\square') = \leg{\lambdamingamma}, \qquad \armm_\xi(\square') = \armtwo{\lambdamingamma}. \]
The conversion of the products to the form on the right side of (1) follows immediately.

For (2), we begin in the same way, 
\begin{align} \left[ \jtwo \right]^* &= \left[\prod\limits_{\square \in \mathrm{C}_{\lambda^-}}\frac{1-q^{\leg{\lambdamingamma}}t^{\armone{\lambdamingamma}+1} }{ 1-q^{\leg{\lambdamingamma}+1}t^{\armone{\lambdamingamma}+1+m_{j-1}(\eta)}}\right]^* \nonumber \\
\label{2j} &= \prod_{\square \in \mathrm{C}_{\lambda^-}} \frac{t^{\armone{\lambdamingamma}+1}-q^{\leg{\lambdamingamma}}}{t^{\armone{\lambdamingamma}+1+m_{j-1}(\eta)}-q^{\leg{\lambdamingamma}+1}}\cdot \prod_{j=1}^{\lambdagamma_{n-k}+1}t^{m_{j-1}(\eta)} \end{align}

The $t$-product counts the columns in $\eta$ which are lower than the newly-added box $x$, which are equivalently those boxes in $\RL$, so the product can be simplified to $t^{|\RL|}$. Consider $\square = (i,j)$ in $\mathrm{C}_{\lambda^-}$ and the corresponding $\square' = (h,j)$ in $C_{\xi+x,\xi}$. Before we compare arm and leg counts, we note that  $m_{j-1}(\eta) = b(\square')$, since the labeled boxes in row $j$ are in the same columns that are nonsymmetric and have height exactly $j-1$ when considered in $\lambdamingamma$. It is straightforward to see that
\[ \legg_{\xi+x}(\square') = \leg{\lambdamingamma} +1,\]
since $\legg_{\xi+x}(\square')$ also counts the box $x$. For the arms,
\[ \armone{\lambdamingamma} = \armm_{\xi+x}(\square')-m_{j-1}(\eta),\]
as both arms count symmetric columns of height $j,\dotsc, \lambdatildegamma_{n-k}$ and nonsymmetric columns of height $j,\dotsc, \lambdamingamma_{n-k}-1$, but only $\armm_{\xi+x}(\square')$ counts nonsymmetric columns of height $j-1$. Thus \eqref{2j} can be converted to
\[t^{|\RL|} \prod_{\square' \in C_{\xi+x,\xi}} \frac{t^{\armm_{\xi+x}(\square)+1-b(\square)} - q^{\legg_{\xi+x}(\square)-1}}{t^{\armm_{\xi+x}(\square)+1} - q^{\legg_{\xi+x}(\square)}}, \]
which completes (2).

For (3), we need to account for the remaining missing monomials. According to \eqref{pIoneComparison},
 $\tildepIone =\left[ p'_{I_1} \right]^* \cdot x \cdot (w_{t_r})^{-1}.$
We also use
  $\left[ \, \overline{\lambdatildegamma}_{n-k}\right]^* = x.$
Then comparing everything,
\begin{align*} \left[ p'_{I_1} \cdot q^{-\mutilde_h+1}\cdot \overline{\lambdatildegamma}_{n-k}\right]^* &= w_{t_r} \cdot q^{-\mutilde_h+1} \cdot \tildepIone\\
&=  t^{c_r}\cdot \tildepIone,\end{align*}
as $-\mutilde_h+1$ is equal to the row index of $w_{t_r}$.\qedhere
\end{proof}

Finally, we observe that in the setting of Theorem~\ref{T:Pieri-match},
\[ c_r = n(\mathrm{sort}(\mu,\eta)) - n(\mathrm{sort}(\lambda,\gamma)),\]
and this completes the proof of part (2) of \cref{T:main}.

\end{document}